\documentclass{article}
\usepackage{amsmath}
\usepackage{amsthm}
\usepackage{amsbsy}
\usepackage{amssymb}
\usepackage{enumerate}
\newtheorem{cor}{Corollary}
\newtheorem{lem}{Lemma}
\newtheorem{thm}{Theorem}
\newtheorem{dfn}{Definition}

\newcommand{\norm}[1]{\lVert#1 \rVert}
\newcommand{\abs}[1]{\lvert#1 \rvert}

\newcommand{\sca}[2]{\ensuremath{\langle #1,#2 \rangle}}
\newcommand{\algebra}[1]{\ensuremath{\mathbb{C}\langle #1 \rangle}}
\newtheorem*{myth}{Theorem}
\DeclareMathOperator{\rad}{R^*}
\begin{document}
\title{On $O^*$- representability and $C\sp*$- representability of  $*$-algebras.}

\author{Stanislav Popovych \\\vspace{0.1cm}\\
{\footnotesize \sl Department of Mathematics, Chalmers University of
Technology, SE-412 96 G\"oteborg, Sweden} \\
{\footnotesize \sl stanislp\symbol{64}math.chalmers.se}}

\maketitle

\footnotetext{ The
author was supported by the Swedish Institute.}

\begin{abstract}
Characterization of the *-subalgebras in the algebra of bounded
operators acting on Hilbert space is presented. Sufficient
conditions for the  existence of a faithful representation in
pre-Hilbert space of a *-algebra in terms of its Groebner basis are
given. These
 conditions are  generalization of the unshrinkability of
 monomial *-algebras introduced by C. Lance and P. Tapper.
 The applications to *-doubles, monomial *-algebras and
several other classes of *-algebras
 are presented.  \end{abstract}

\section{Introduction}\label{intro}

 In the paper we study  conditions for  a  $*$-algebra
to be faithfully represented by bounded or unbounded operators on a Hilbert
space.

The term "algebra of unbounded operators" admits different
interpretations. In present work this term means $O\sp*$-algebra
(\cite[p.36]{smudgen}), i.e.  a $*$-subalgebra of the algebra of
linear operators acting on a pre-Hilbert space.  Let  $\rm E$ denote
a pre-Hilbert space and  $\rm H$ a Hilbert space which is the
completion of $\rm E$. The $*$-algebras of linear operators acting
on these spaces are denoted by $\rm L(E)$ and $\rm L(H)$. Let $A$ be
a $*$-algebra over complex numbers. In this paper  we study
conditions for the existence of an embedding of $A$ into $\rm L(E)$
and $\rm L(H)$. In the first case, it is equivalent to $A$ being
$*$-isomorphic to a $O\sp*$-algebra, such algebras will be called
$O\sp*$-representable. In the second case $A$ is isomorphic to a
pre-$C^*$-algebra and we will say (following C. Lance and
P.Tapper~\cite{lance}) that $A$ is $C^*$-representable.

If  $A$ is embedded in $\rm L(E)$ and every operator
$a\in A$ is bounded then one can extend each  $a\in A$ to an
operator acting on  $H$  and thus obtain an
inclusion $A\hookrightarrow \rm L(H)$. In the  general case $A$ will
be represented by unbounded   operators  on $H$ such that the
intersection of their  domains is dense.

The celebrated Gelfand-Naimark theorem characterizes closed $*$-subalgebras of $\rm L(H)$ in terms of the norm on a $*$-algebra. There are also characterizations of such subalgebras in terms of orders on the set of self-adjoint elements~\cite{Ajupov}. A related result was obtained by  J. Kelly and R. Vaught in \cite{Kelley}. They describe  universal $C\sp*$-seminorm for Banach $*$-algebras ${A}$ with isometric involution.

 The non-complete subalgebras of $\rm L(H)$  are less well  studied. A characterization of pre-$C\sp*$-algebras inside  the class of normed $*$-algebras is given by G. Allan (see~\cite[p. 281]{DoranBelfi}). The algebraic characterization of $C\sp*$-representability of T$^*$-algebras was obtained in \cite{pop1} and in \cite{Palmer} by different methods. In the next section we derive this  result as a simple consequence of an abstract characterization of matrix ordered $*$-algebras recently obtained  in \cite{juspop}. The latter is analogues  to the Choi and Effors characterization  of abstract operator systems.
In section \ref{scbound} we present sufficient conditions of operator representability in algebraic terms. In particular, we seek a generalization of a particularly simple necessary condition  of the $C\sp*$-representability of a $*$-algebra $A$ that the equation $x^*x=0$ has only zero solution in $M_n(A)$ for all $n\ge 1$. Algebras possessing this property are called \textit{completely positive} in the present paper and \textit{ordered} in \cite{Palmer}.  We  prove (see Corollaries~\ref{simple} and \ref{directsum}) that for a large class of $*$-algebras complete positivity is also sufficient for  $C\sp*$-representability.  We also present   example using Gr\"obner bases theory which shows that the condition of complete positivity is not sufficient in general.

  The literature on the $O\sp*$-representability of
  finitely presented $*$-algebras consists so far  of
  isolated classes of examples.
  In~\cite{pop2}, the author proved that a monomial $*$-algebra is
  $O^*$-representable if and only if in the minimal defining set of
  monomial relations of  the form $w_j=0$ where $w_j$ is a word,
  all $w_j$ are unshrinkable. It should be noted that
   Lance and Tapper~\cite{tapper,lance} conjectured that
   such $*$-algebras are  $C^*$-representable.
   This is still an open problem.
In Section~\ref{class} we introduce a larger
class of $O\sp*$-representable  $*$-algebras which we call
non-expanding (see Definition~\ref{nonexpand}).
This class is a generalization of monomial $*$-algebras.
The main novelty of our approach is that we use the notion of
Gr\"obner basis  to define this class and use methods of
Gr\"obner bases theory to establish $O\sp*$-representability
and derive further results.

The sufficient conditions of non-expandability obtained
in Section~\ref{suffic} allowed one to  show
that several known classes of $*$-algebras
fall in the class of non-expanding  $*$-algebras.
Thus their representability could be treated from a
unified point of view. These sufficient conditions
are algorithmically verifiable for $*$-algebras
given by a finite number of generators and relations.

\section{  Representability by bounded operators.}\label{scbound}
In this section  several characterizations  of representability of a $*$-algebra  by bounded operators acting on a Hilbert space $\rm H$ are presented. These characterizations are consequences of the results of \cite{juspop}.

 If a $*$-algebra $A$ is $*$-isomorphic to a subalgebra of a $C^*$-algebra $\mathcal{A}$ then by the Gelfand-Naimark theorem $A$ is also $*$-isomorphic to a subalgebra of $\rm L(H)$ and thus can be faithfully represented by bounded operators on $\rm H$. Such  $*$-algebra is called $C^*$-representable (see~\cite{lance}).

Let  ${A}_{sa}$ denote the set of self-adjoint elements in $A$.
The following definition was introduced in~\cite{powers}.
\begin{dfn} Given a $*$-algebra $A$ with unit $e$, we say that a subset $C\subset A_{sa}$ is
\textit{algebraically admissible cone} if
\begin{enumerate}
\item $C$ is a cone  in $A_{sa}$ and  $e\in C$; \item $C \cap (-C) =\{0\}$;
\item $x C x^* \subseteq C$ for every $x\in A$;
\end{enumerate}
\end{dfn}

With a cone $C$ we can associate a partial order $\ge_C$ on the real vector space $A_{sa}$ given by the rule  $a\ge_C b$ if $a-b\in C$.  Henceforth we will suppress subscript $C$ if it  will not lead to ambiguity.

\begin{dfn}
 Recall that an element $u \in A_{sa}$ is called an \textit{order unit} for $A_{sa}$ provided that for every $x\in A_{sa}$ there exists a positive real $r$ such that $r u + x \in C$. An order unit $u$ is called \textit{Archimedean} if $r u + x\in C$ for all $r>0$ implies that $x\in C$.
\end{dfn}

\begin{myth}[cf. \cite{juspop}]
A $*$-algebra $A$ with unit $e$ is $C^*$-representable if and only if 
there is an algebraically admissible cone on $A$ such that $e$ is an Archimedean order unit.
\end{myth}

The assumptions of the $C^*$-representability criterion given
in the above theorem  are the same as  in  Choi and Effros
characterization of abstract operator systems~\cite{ChoiEffros},
however no additional structure on the matrices is required
and the matrix order is replaced with the order given by an algebraically admissible cone.

For  $*$-algebra $A$  and  algebraically admissible cone $C$ with order unit $e$ the function
$\norm{a} =  \inf \{r>0: r e \pm a \in C \}$ is a seminorm on $A_{sa}$
(see \cite[lemma 4]{juspop}). The function  $\abs{x} = \sqrt{\norm{x^*x}}$ is a $C\sp*$-seminorm on $A$ and for self-adjoint $a\in A$ we have $\norm{a}\le \abs{a}$ (see \cite[lemma 5, Theorem 6]{juspop}).

The main drawback of the characterization given  in the above theorem is that it requires some additional structure on a $*$-algebra. Our next  objective is to give an intrinsic characterization of
$C^*$-representability using the algebraic structure of $*$- algebra alone. It turns out
to be possible in case the $*$-algebra is bounded. For such $*$-algebras the set of positive elements forms an  algebraically admissible cone.

\begin{dfn}
Recall that a $*$-algebra $A$ is called \textit{$*$-bounded} if  for every $a\in A$ there is constant $C_a$  such that for every  $*$-representation $\pi: A\to B(H)$ we have $\norm{\pi(a)}\le C_a$.
\end{dfn}

\begin{dfn}
 An element $a \in {A}_{sa}$  is called \textit{positive} if  $a=\sum_{i=1}^{n}a^{*}_{i}a_{i}$ for some $n\ge 1$ and $a_i\in A$ for $1\le i\le n$. The set of positive  elements in $A$ will be denoted by ${{A}}_{+}$.
\end{dfn}

It is easy to check that  the  cone $A_+$ of a unital $*$-algebra $A$ is  algebraically  admissible. To formulate our next result we will need some definitions from the theory of ordered algebras (\cite{Ajupov}).
\begin{dfn}
Let $A$ be a unital $*$-algebra.
\begin{enumerate}
\item An element  $a \in {A}_{sa}$ is \textit{bounded} if there is
$\alpha \in \mathbb{R}_+$ such that $\alpha e\ge a \ge -\alpha e$.
\item An element  $x=a+ib$ with $a,\ b \in {A}_{sa}$ is \textit{bounded} if so are the elements $a$ and  $b$.
\item   The algebra ${A}$ is \textit{bounded} if all its elements are bounded.
\end{enumerate}
\end{dfn}

We will collect some useful facts about bounded elements in the following Lemma.  They can be found in ~\cite[proposition 1, p. 196]{Ajupov}:
\begin{lem}\label{boundelem}
Let $A$ be a unital $*$-algebra then
\begin{enumerate}
\item the set of all bounded elements in  ${A}$ is a  $*$-subalgebra in $A$;
\item an element  $x\in {A}$ is bounded if and only if  $xx^*$ is bounded;
\item if  ${A}$  is  generated by a set $ \{ s_j \}_{j\in J} $ such that each  $s_js_j^*$ is  bounded then  ${A}$ is bounded.
\end{enumerate}
\end{lem}

For example, an  algebra  ${A}$  generated by  isometries (i.e., elements satisfying relation  $s^*s=e$) or projections  (i.e., elements satisfying relation $p^*=p=p^2$) is bounded.
 One can easily prove that a bounded  $*$-algebra ${A}$ is $*$-bounded and thus there exists its  universal enveloping $C^*$-algebra $C\sp*(A)$.

Recall the definition of $*$-radical introduced by Gelfand and Naimark (see for example \cite[(30.1)]{DoranBelfi}).
\begin{dfn}
For a $*$-algebra $A$ the $*$-\textit{radical} is the set   $\rad(A)$ which is the intersection of  the  kernels of all topologically irreducible  $*$-representations of  $A$ by bounded operators on Hilbert spaces.
\end{dfn}

It is  known that $\rad(A)$ is  equal to the intersection of  the kernels of all $*$-representations (see for example~\cite[Theorem (30.3)]{DoranBelfi}). Clearly the factor algebra ${A}/\rad({A})$ of a $*$-bounded algebra $A$ is $C^*$-re\-pre\-sent\-able.

 The following theorem provides an intrinsic characterization of $C\sp*$-representability of
bounded $*$-algebras.

\begin{thm}\label{raddesc}
Let  ${A}$ be a bounded $*$-algebra then the following holds.
\begin{enumerate}
\item $\abs{x}$ coincides with the norm of the universal enveloping $C^*$-algebra $C\sp*(A)$ of $x\in {A}$,  i.e. $\abs{x}=sup_{\pi} \norm{\pi(x)} $ where  $\pi$ runs over all $*$-representations  of ${A}$ by bounded operators on Hilbert spaces. Thus
\[
\sup_{\pi} \norm{\pi(x)}^2=\inf_{f\in  {A}_+}\{(xx^*+f)\cap\mathbb{R}e\}.
\]
Moreover,  $\norm{a} = \abs{ a}$ for self-adjoint $a\in A$.
\item The null-space of $|\cdot|$ which is  $\rad({A})$  consists of those  $x$ such that for every   $\varepsilon > 0$ there are  $x_1, \ldots, x_n$ in  ${A}$ satisfying the equality \begin{equation}\label{rad}
x^*x+ \sum_{j=1}^{n}x_j^*x_j=\varepsilon e.
\end{equation} \item  $A$ is $C\sp*$-representable if and only if  $\rad{(A)}=\{0\}$.
\end{enumerate}
\end{thm}

\begin{proof}

Since every  $x$ in $A$ is bounded there are real  $\alpha>0$ and  $x_1,\ldots, x_m$  in $A$ such that
\begin{equation}\label{eq1}
xx^*+\sum_{i=1}^m x_ix_i^*=\alpha.
\end{equation}
If  $\pi$ is a representation of
${A}$ by bounded operators then  $\norm{ \pi(xx^*)} \le \alpha$. Thus
$\sup_\pi \norm{ \pi(x) }^2\le \inf \alpha$, where  $\pi$ runs over all $*$-representations of ${A}$ and infimum  is taken over all $\alpha$ as in (\ref{eq1}). Therefore   $\abs{x}\ge \sup_{\pi} \norm{\pi(x)}$ for all $x\in A$. The converse inequality  also holds since the right-hand side is the maximal pre-$C^*$-norm.  This proves the universal property of the pre-norm  $|\cdot|$. Obviously its null-space is $\rad({A})$. By \cite[Lemma~5]{juspop}, $\norm{a}\le \abs{a}$ for every self-adjoint $a\in A$. But inequality $-\alpha e \le a \le \alpha e$ implies that $-\alpha {\rm I} \le \pi(a) \le \alpha{\rm I}$ for every $*$-representation $\pi$ and identity operator ${\rm I}$. Hence $\norm{\pi(a)} \le \alpha$. From this follows  $\abs{a}\le \norm{a}$ and, consequently, $\abs{a} =  \norm{a}$.

  Thus we only have to prove that the null-space of $|\cdot|$ is the set of all $x$ such that for every $\varepsilon>0$ there are $x_1, \ldots, x_n$ in $A$ such that (\ref{rad}) is fulfilled. The null-space is the set of $x$ such that $\inf \{ r>0: re \ge x^* x \ge -r e\} =0$. But by definition of the order  $r e -  x^*x \ge 0$ if there $x_1,\ldots, x_n \in A$ such that $r e -  x^*x  = x_1^*x_1+\ldots x_n^* x_n$ which proves $(2)$ and the theorem. \end{proof}

Note that J. Kelly and R. Vaught in 1953 proved that
\begin{equation}\label{st}
\sup \norm{\pi(x)} = \inf\{t\in \mathbb{R}_+| t^2 -x^*x \in {A}_+\}
\end{equation} where  ${A}_+ = \{\sum_{j=1}^n a_j^*a_j, n\in\mathbb{N}, a_j\in {A}\}$,
 $\pi$ runs over all
$*$-representations  for Banach $*$-algebras ${A}$ with isometric
involution (see~\cite{Kelley}).  The proof of formula \eqref{st} based on the Hahn-Banach theorem
  for any T$^*$-algebra (a $*$-algebra $A$ such that every $x\in A_{sa}$ is bounded) presented in monograph~\cite{Palmer} and by purely algebraic methods in \cite{pop1}.

As a corollary of the above theorem we  obtain the following description of the elements positive in every representation.
\begin{cor}
Let $A$ be a bounded $*$-algebra. An element $a\in A_{sa}$ has the property that $\pi(a)\ge 0$ for each $*$-representation $\pi$ of $A$ in $\rm L(H)$ if and only if for every $\varepsilon>0$ there are $x_1, \ldots, x_n\in A$ such that  $a+\varepsilon = \sum_{j=1}^n x_jx_j^*$.
\end{cor}
\begin{proof} Clearly, given $a\in A$, $\tau(a)\ge 0$ for every $*$-representation $\tau$ of $A$ in $\rm L(H)$ if and only if  $\pi(a)\ge 0$ for universal representation $\pi$ of $A$. Since every  representation could be factored through the universal representation  $\pi$,   $\abs{x} = \norm{\pi(x)}$ for all $x\in A$. Here $\abs{\cdot}$ is the  norm  as in  Theorem~\ref{raddesc}.    A self-adjoint operator $\pi(a)$ is positive if and only if $\norm{ C {\rm I} -a }\le C$ where $C= \norm{\pi(a)}$ and ${\rm I}$ is the identity operator. Thus assuming  $\pi(a)\ge 0$  we have, by Theorem~\ref{raddesc}, that $\abs{ \abs{a} -a } \le \abs{a}$ and hence  $\norm{ \abs{a} -a } \le \abs{a}$. Consequently, $ \abs{a} -a \le \abs{a}+\varepsilon$  for every $\varepsilon>0$. Which means that $a+\varepsilon$ can be written as $\sum_{j=1}^n x_jx_j^*$ for some $x_j\in A$. The converse statement is obvious.
\end{proof}

It is  well known  that for a finite dimensional $*$-algebra $A$ the necessary and sufficient conditions for $C^*$-representability is that $A$ is \textit{positive}, i.e. the equation $x^*x = 0$ has only zero solution in $A$. For an infinite dimensional $*$-algebra $A$ the above condition is not sufficient  since there are positive (even commutative) $*$-algebras such that $M_2(A)$ is not positive (see \cite[Example (32.6)]{DoranBelfi}). This motivates the following definition.

\begin{dfn}
A $*$-algebra $A$ is called \textit{completely positive} if $M_n(A)$ is positive for every $n\ge 1$.
\end{dfn}
We will prove below that for a large class of $*$-algebras the complete positivity is equivalent to $C^*$-representability. However,  we will also present examples of completely positive algebras which are not $C^*$-representable.

Consider the  inductive limit $M_\infty(\mathbb{C}) = \lim (M_n(\mathbb{C}), \phi_n)$  where  $$\phi_n(a)=\left(%
\begin{array}{cc}
  a & 0 \\
  0 & 0 \\
\end{array}%
\right)$$ is an embedding of $M_n(\mathbb{C})$ into $M_{n+1}(\mathbb{C})$.
It is clear that $A$ is completely positive if and only if $A\otimes M_\infty(\mathbb{C})$ is positive. Since the $*$-algebra $ M_\infty(\mathbb{C})$ is not unital and is not finitely generated  we prefer to replace it with the  Teoplitz $*$-algebra $\mathcal{T} = \mathbb{C}\langle  u, u^* | u^* u = e \rangle$ in the above characterization of complete positivity.

\begin{thm}\label{tensor}
For a $*$-algebra $A$ the following conditions are equivalent.
\begin{enumerate}
\item
$A$ is completely positive.
\item For every  $n\ge 1$ the equation $x_1^*x_1+\ldots x_n^* x_n = 0$ has only zero solution $x_1=\ldots =x_n=0$ in $A$. \label{ord}
\item $A\otimes \mathcal{T}$ is positive.
\end{enumerate}
\end{thm}
\begin{proof}
If $x_1^*x_1+\ldots x_n^* x_n = 0$ for some $x_1, \ldots, x_n \in A$ then for a matrix $C\in M_n(A)$ with the first row equal to $(x_1, \ldots, x_n)$ and the rest rows being zero we have $C C^* = 0$. Thus (1) implies (2). If for some non-zero matrix $D\in M_n(A)$ we have $DD^*=0$ and j-row is not-zero then considering $(j,j)$-entry in $DD^*$ we have  $d_{j1}d_{j1}^*+\ldots + d_{j1}d_{j1}^* =0$. Thus (2) is equivalent to (1).

It is easy to see that the element $p=e-u u^*$ is a projection in $\mathcal{T}$ and the elements $e_{ij} = u^{i-1} p (u^*)^{j-1}$ for $i, j\le n$ satisfy  the matrix units relations and thus generate an algebra isomorphic to $M_n(\mathbb C)$. From this it follows that $A\otimes \mathcal{T}$ contains a subalgebra isomorphic to $A\otimes M_\infty(\mathbb{C})$. Hence the condition that $A\otimes \mathcal{T}$ is positive implies that $A$ is completely positive.

We prove now the converse statement. Assume that $A$ is completely positive. Since the relation $u^*u-e$ constitutes a Gr\"obner basis for $\mathcal{T}$  the set $\{u^k u^{*l} | k\ge 0, l\ge 0\}$ forms a linear basis for $\mathcal{T}$. Thus arbitrary $x\in A\otimes \mathcal{T}$ can be written in the form $\sum_{i=1,j=1}^n a_{i,j}\otimes u^i u^{*j}$, where $a_{i,j}\in A$.  Using the  relation  $u^*u=e$ we obtain
\begin{align*}
x^*x = \sum_{i\le k} a_{i,j}^* a_{k,l}\otimes u^j u^{k-i} u^{*l} + \sum_{i'>k'} a_{i',j'}^* a_{k',l'}\otimes u^{j'} u^{k'-i'} u^{*l'}=\\=
\sum_{s=1}^n \sum_{l=1}^n \left[ \sum_{j=1}^s \sum_{k=s-j+1}^n  a_{j+k-s,j}^* a_{k,l} + \sum_{r=1}^l \sum_{i=l-r+1}^n a_{i,s}^* a_{i+r-l,r} \right] u^s u^{*l}.
\end{align*}
 Thus $x^*x=0$ would imply that for every $1\le s, l \le n$:
 \begin{equation}
\sum_{j=1}^s \sum_{k=s-j+1}^n  a_{j+k-s,j}^* a_{k,l} + \sum_{r=1}^l \sum_{i=l-r+1}^n a_{i,s}^* a_{i+r-l,r} = 0. \label{teopl}
\end{equation}
For $s=1$ and $l=1$ we have $\sum_{k=1}^n a_{k,1}^*a_{k,1}+\sum_{i=1}^n a_{i,1}^*a_{i,1}=0$. Since $A$ is completely positive we have $a_{k,1} = 0$ for all $1\le k \le n$.
We will  prove that $a_{k,t}=0$ for all $k$ using  an induction on $t$. We have already check the base of the induction. So assume that $a_{k,m}=0$ for all $k$ and prove that $a_{k,m+1}=0$. Setting $s=l=m+1$ in (\ref{teopl}) and using the induction hypothesis we obtain
\begin{align*}
\sum_{k=1}^n  a_{k,m+1}^* a_{k,m+1} + \sum_{i=1}^n  a_{i,m+1}^* a_{i,m+1} = 0.
\end{align*}
Since $A$ is completely positive we get $a_{k,m+1}=0$ for all $1\le k \le n$ which proves our induction claim and the theorem.
\end{proof}

Note that $*$-algebras satisfying conditions (\ref{ord}) are called \textit{ordered} in \cite{Palmer}.
We prefer the term completely positive to avoid confusion with the order given by cones and to emphasize the analogy with completely positive maps.

One can easily show   that   complete positivity  is preserved under taking sub-direct products,    direct limits and taking subalgebras. It also preserved under making extensions, i.e. if $J$ is a $*$-ideal in $A$ which, considered as $*$-algebra,  is completely positive and such that $A/J$ is also completely positive then $A$ itself is completely positive. Indeed, if $\sum_{j=1}^n  x_j^*x_j= 0$ in $A$ then passing to the factor algebra $A/J$ and using its complete positivity  we obtain that $x_j$ are elements from $J$. Using completely positivity  of $J$ we conclude that $x_j=0$ for all $j$.

It is an open question  whether the tensor product $A\otimes B$ of two completely positive $*$-algebras is completely positive.  However, it can be easily checked that a tensor product of two of two $O\sp*$-representable algebras is $O\sp*$-representable.

Using Theorem~\ref{tensor} we can simplify the conditions of  Theorem~\ref{raddesc} in the following way.

\begin{thm}\label{useteopl}
Let $A$ be a bounded unital $*$-algebra and $\mathcal{T}$ be the Teoplitz $*$-algebra. Then $A$ is $C\sp*$-representable if and only if every     $x\in A\otimes \mathcal{T}$ with the property that   for every $\varepsilon>0$ there exists $y\in A\otimes \mathcal{T}$ such that $xx^* +yy^* = \varepsilon (e- uu^*)$ is  zero.
\end{thm}
\begin{proof}

To prove that $A$ is $C\sp*$-representable it is suffices to prove that $\rad{(A)} = \{0\}$.
 If $x\in \rad{(A)}$ then, by  Theorem~\ref{raddesc},  for every $\varepsilon>0$ there are $x_1, \ldots$, $x_n \in A$ such that $xx^* + \sum_{j=1}^n x_j x_j^* = \varepsilon e$. Consider   $n\times n$-matrices $X$ and $C$ with coefficients in  $A$ such that the first row of  $X$ is   $(x,0,\ldots, 0)$  and the first row of  $C$ is $(x_1,x_2, \ldots,x_n)$  and all other rows of $X$ and $C$ are equal to zero.

Since the  subalgebra $B_n$ of $\mathcal{T}$ with basis $e_{ij}$ is isomorphic to $M_n(\mathbb{C})$. One can identify  $B_n$ with $M_n(\mathbb{C})$ and consider the  algebra  $M_n(A)\simeq A\otimes M_n(\mathbb{C})$ as a subalgebra of $A\otimes \mathcal{T}$. Moreover, after this identification one has $XX^* + CC^* = \varepsilon (e-uu^*)$. Thus $X=0$ and, consequently, $x=0$.

The necessity of the conditions of the theorem follows easily from the fact that $\mathcal{T}$ is $C^*$-representable and thus its tensor product with any $C^*$-representable algebra $A$ is also  $C^*$-representable.
\end{proof}
\begin{cor}
Each bounded completely positive  $*$-algebra $A$ has a non-trivial representation in $B(H)$.
\end{cor}
\begin{proof} Assume that  $\abs{e}=0$. Then there are $x_1, \ldots, x_m \in {A}$ such that $e+x_1x_1^*+ \ldots +x_mx_m^*=\frac{1}{2}e$. Therefore $\sum_{j=1}^{m} x_jx_j^* +yy^*=0$ where $y=\frac{1}{\sqrt{2}}e$, which contradicts the complete positivity  of $A$.
Hence  $\abs{e} \ne 0$. For  the universal representation  $\pi$ of  $A$, which is a faithful representation of the enveloping $C^*$-algebra $C\sp*({A})$, we  have  $\pi(e) \ne 0$.
\end{proof}
The assumptions of the previous corollary  can be  weakened.
Recall that  an ideal  $I$ of a $*$-algebra $A$ is called {\it endomorphically closed} if  $f(I)\subseteq I$ for every $*$-endomorphism  $f:A\to A$. An algebra
 $A$ is called {\it endomorphically simple} if it has only trivial endomorphically closed $*$-ideals. We will say that a $*$-ideal $J$ of $A$ is {\it square root  closed} if for every elements   $x_1,\ldots, x_n \in A$ equality $\sum_{j=1}^n x_j x_j^*\in J$  implies that $x_j \in J$. This is equivalent to $A/J$ being completely positive.

\begin{cor}\label{simple}
Let  $A$ be a bounded  unital $*$-algebra without non-trivial  endomorphically closed and square root  closed $*$-ideals. Then
 $A$ is   $C^*$-representable if and only if
${A}$ is completely positive.
\end{cor}
\begin{proof} The necessity is obvious. Since  the $*$-radical of  a  $*$-algebra is an endomorphically closed and a square root  closed $*$-ideal which, by the previous corollary, does not coincide with  $A$, it must be zero.
\end{proof}
\begin{cor}\label{directsum}
If a unital bounded algebra  ${A}$ is a direct sum of endomorphically simple $*$-algebras ${A}_n$, then  ${A}$ is  $C\sp*$-re\-pre\-sent\-able if and only if $A$ is completely positive.
\end{cor}
\begin{proof}
Let $\pi_n$ be the canonical  $*$-homomorphism   ${A}\to {A}_n$. By Lemma~\ref{boundelem}, for any $a\in A$, there are elements $a_j \in A$ and $c \in \mathbb{R}$  such that  $ c e - a^*a = \sum_{i=1}^{n}a^{*}_{i}a_{i}$. Thus $c e - \pi_n(a) \pi_n(a)^*$ is a positive element of $A_n$. Hence $\abs{\pi_n(a)\pi_n(a)^*} < c$. Since $\pi_n$ is subjective $A_n$ is bounded by Lemma~\ref{boundelem}.  The previous corollary then imply that each ${A}_n$ is  $C^*$-representable and hence the same is true for their direct sum ${A}$.
\end{proof}

\begin{thm}\label{twofunc}
A bounded  $*$-algebra $A$ is $C^*$-representable if and only if there  are  mappings $F:A_{+}\to \mathbb{R}$ and $G:A_{+}\to \mathbb{R}$ such that

1. $F(aa^*)>0$ for each   $a\ne 0$

2. $G(\sum_{i=1}^{n}a_ia_i^*)\ge F(a_ja_j^*)$ for arbitrary elements    $a_1, \ldots, a_n \in A$ and every  $j\in\{1,\ldots,n\}$.

3. $\lim_{\varepsilon\to 0+}G(\varepsilon e)=0$ for $\varepsilon\in\mathbb{R}$.
\end{thm}

\begin{proof}
If $A$ is not  $C^*$-representable, then there is a nonzero  $x\in \rad(A)$. By Theorem~\ref{main}, for each
$\varepsilon >0 $ one can find  $x_1,\ldots, x_l \in A$ such that   $xx^*+\sum^l_{i=1}x_ix_i^*=\varepsilon e$ and thus
  $G(\varepsilon e)\ge F(xx^*)$. From this we obtain
 $F(xx^*)=\lim_{\varepsilon\to 0} G(\varepsilon e)=0$ contrary to   the condition~1 of the theorem.

If  $A$ is  $C^*$-representable then there is pre-$C^*$-norm   $\norm{\cdot}$ on $A$. Put  $G(x)=F(x)=\norm{x}$. For each positive  $x$ in $A$,  $F(x)=\sup s(x)$ where supremum is taken over all states on the enveloping $C^*$-algebra  $C^*(A)$. For every state  $s$ we  have $s(\sum_{i}x_ix_i^*)\ge s(x_jx_j^*)$ and, taking supremum, we obtain   $G(\sum_{i}x_ix_i^*)\ge F(x_j x_j^*)$
\end{proof}

Now we  apply  Theorem~\ref{raddesc} to the group $*$-algebras. Let $G$ be a discrete group and  $\mathbb{C}[G]$  its group $*$-algebra. Elements of $\mathbb{C}[G]$ could be considered both as a formal linear combinations of elements of $G$ with complex coefficients and as a functions from $G$ to $\mathbb{C}$ with finite support. Let $P$ denote the set $\{\sum_{j=1}^n f_jf_j^*|n\in \mathbb{N},  f_j\in \mathbb{C}[G] \}$ which is a subset of the set of positive definite functions on $G$ with compact support. By important result due to Godement \cite[(13.8.6)]{Dix} each element of $P$ is of the form $f f^*$ for some $f\in \mathbb{C}[G]$.

Considered as a positive definite function element $\phi\in P$ give rise to a cyclic representation $\pi_\phi$ in a Hilbert space with cyclic vector $\xi$ such that $\phi(s) = (\pi_\phi (s) \xi,\xi)$ for every $s\in G$. By~\cite[Lemma 14.1.1]{Dix} for every $f\in \mathbb{C}[G]$ and $\phi \in P$ we have that $\norm{\pi_\phi(f)} \le \norm{\lambda(f)}$ where $\lambda$ denote left regular representation of $\mathbb{C}[G]$. Since $\delta_e\in P$ and $\pi_{\delta_e}= \lambda$,  $\sup
 _{\phi\in P}\norm{\pi_\phi(f)} =  \norm{\lambda(f)}$. Thus using the set $P$ one can recover  the norm of the reduced group $C^*$-algebra $C^*_{red} (G)$. By the next corollary it is also possible to recover  the norm of the  group $C^*$-algebra $C^*(G)$.
\begin{cor}
  Let $\norm{\cdot}$ denote the norm on $C^*(G)$. Then for every $f\in \mathbb{C}[G]$ the following formula holds
$$\norm{f}^2= \inf_{\phi\in P}\{(\phi +ff^* )\cap\mathbb{R}e\}.$$
\end{cor}
\begin{proof}
 Clearly $P$ is the set of positive elements of $*$-algebra $\mathbb{C}[G]$. For every $f\in \mathbb{C}[G]$ norm $\norm{f}$ is the norm of universal enveloping $C^*$-algebra of $\mathbb{C}[G]$ and consequently, by Theorem~\ref{raddesc}, $\norm{f}^2= \inf_{\phi\in P}\{(\phi +ff^* )\cap\mathbb{R}e\}$. \end{proof}

Since $G$ is amenable if an only if for every $f\in \mathbb{C}[G]$ reduced norm is equal to universal enveloping norm  we obtain the following.
\begin{cor} A discrete group $G$ is non-amenable if an only if there exists $f\in \mathbb{C}[G]$ and $\varepsilon>0$ such that for every $g\in\mathbb{C}[G]$ element $\frac{\norm{fg}_2}{\norm{g}_2} +\varepsilon $ can not be presented in the form $ff^* +  g g^*$ for some $g\in \mathbb{C}[G]$. Here $\norm{g}_2^2= \sum_{k=1}^m \abs{\alpha_k}^2$ for the element $g= \sum_{k=1}^m \alpha_k w_k$ with $\alpha_k\in\mathbb{C}$ and distinct $w_k\in G$.
\end{cor}
 In the following  example we present  a completely positive  bounded $*$-algebra which is not  $C^*$-representable.
For the definitions of the Gr\"obner basis,  the set of basis words $BW$ and the operator ${\rm R_S}$  we refer the reader to the appendix.

{\bf Example. } Consider $*$-algebra given by generators and relations $$A=\mathbb{C}\big<a,\ x|a^*a=qaa^*, xx^*+aa^*=e\big>$$ where parameter $0<q<1$.
Clearly, $A$ is bounded. It can be easily checked that the set $S=\{a^*a - qaa^*, xx^*- aa^*- e\}$ is a Gr\"obner basis of $A$. Thus the set $BW$ consisting of the words containing no subword from the set
 $\{a^*a, xx^*\}$ forms a linear basis for $A$.
For arbitrary   $z$ in $\mathbb{C}\big<a,\ x\big>$ the element ${\rm R_S}(z)$ could be written as  $\sum_{i=1}^n \alpha_iu_ix^{k_i}$, where $u_i$ does not end with  $x$,  $k_i \ge 0$, $\alpha_i\not= 0$ and $u_i\in BW$ for all $1\le i\le n$.

Let $t$ be the minimal length of the words $u_ix^{k_i}$. Put  $J=\{j: |u_j| = t \}$. Denote by   $F(z)$ the sum of those  $\alpha_i$ with  $i\in J$  such that    $u_ix^{k_i}=ww^*$ for some word $w$.
 We will prove   that
$F(zz^*)=\sum_{j\in J}|\alpha_j|^2$. Indeed,   \begin{gather*}{ \rm R_S}(u_ix^{k_i}x^{*k_j}u_j^*) = \\     \begin{cases}
-u_i(\sum_{1\le s\le k_i}x^{k_i-s}aa^*x^{*k_j-s})u_j^* + {\rm R_S}(u_iu_j^*), &\text{ if $k_i=k_j$} \\
-u_i(\sum_{1\le s\le min(k_i,k_j)}x^{k_i-s}aa^*x^{*k_j-s})u_j^*,
& \text{ if  $k_i\not=k_j$ }
\end{cases}
 \end{gather*}
The sum $u_i(\sum_{1\le s\le min(k_i,k_j)}x^{k_i-s}aa^*x^{*k_j-s})u_j^*$ contains no words of length $t$. Thus computing $F(zz^*)$ it is sufficient to consider only the sum
$$-u_i(\sum_{1\le s\le k_i}x^{k_i-s}aa^*x^{*k_j-s})u_j^* + {\rm R_S}(u_iu_j^*).$$
Since both  $u_i$ and  $u_j$ do not end with
 $x$ the element ${\rm R_S}(u_iu_j^*)$ is a monomial of length  $\abs{u_i}+\abs{u_j}$. Thus, if some monomial ${\rm R_S}(u_iu_j^*)$  in  ${\rm R_S}(zz^*)$ has minimal length (which is equal to $2t$)  then  $i,\ j\in J$ (in particular  $\abs{u_i}=\abs{u_j}$). Equality   $u_iu_j^*=ww^*$ implies $u_i=u_j$. Indeed, if
 $u_i$ ends with  $a$ or with $x^*$ or word $u_j$ ends with $a^*$
or with $x^*$ then  ${\rm R_S}(u_iu_j^*)$ is just $u_iu_j^*$ (as in free $*$-algebra). Thus using  equality $u_iu_j^*=ww^*$ we can conclude that $u_i=u_j$. Otherwise, write  $u_i=v_ia^{*k}$ and
$u_j=v_ja^{m}$ where $v_i$ does not end with  $a^*$ and   $v_j$ does not end with  $a$. Thus  ${\rm R_S}(u_iu_j^*)=q^{km}v_ia^ma^{*k}v_j^*$.
If $m>k$ then, since  $u_iu_j^*=ww^*$, we have
$
v_ia^{m_1}=w$ and $a^{m-m_1}a^{*k}v_j^*=w^*$, for some $1\le m_1<m$.
This is a contradiction since  $w$ ends with $a$ and  $a^*$ simultaneously.   Similarly if  $m<k$ then
$
w=v_ia^ma^{*k_1}$ and  $w^*=a^{*(k-k_1)}v_j^*$, for some  ($1\le k_1<k$). We obtain that $w$  ends with $a$ and $a^*$ which is again  a contradiction. Thus  $m=k$ and  $w=v_ia^{k}=v_ja^{k}$. So
$v_i=v_j$ and $u_i=u_j$. We have proved  so far that
  $u_iu_j^*=ww^*$ implies that $u_i=u_j$. From this it easily follows that $F(zz^*)=\sum_{j\in J}|\alpha_j|^2$.
Obviously    $F(aa^*)>0$ if  $a\ne 0$ and
\[
F(\sum_{i=1}^n a_ia_i^*)\ge \min_{i}F(a_ia_i^*),
\] end clearly $ F(\varepsilon e)=\varepsilon$ for $\varepsilon \in \mathbb{R}$.
Thus  $A$  is  completely positive  $*$-algebra. If $\pi$ is a representation of $A$ in Hilbert space then
\[
\norm{\pi(aa^*)}=\norm{\pi(a^*a)}=q\norm{\pi(aa^*)},
\] which implies that
$\norm{\pi(aa^*)}=0$. Thus  $A$ is not $C^*$-representable.

\section{ Unshrinkable words and  Gr\"obner bases. }\label{class}
C. Lance and P. Tapper (cf.~\cite{lance, tapper})  studied $C^*$-representability of
$*$-algebras $A_w$  generated by  $x$ and $x^*$ with one monomial defining
relation $w=0$ where $w=x^{\alpha_1} x^{*\beta_1}\ldots x^{\alpha_k} x^{*\beta_k}$, $\alpha_j$ and $\beta_j$ are positive integers. They conjectured that $A_w$ is $C^*$-representable if and only if the word $w$ is  \textit{unshrinkable}, i.e. $w$ can not be presented in the form $d^*d u$ or $ud^*d$ where $u$ and $d$ are words and $d$ is non-empty. A very appealing feature of this conjecture is that being true it gives a condition of $C^*$-representability of a monomial $*$-algebras in terms of its defining relations.   In~\cite{pop2} the author proved that a  monomial $*$-algebra is $O\sp*$-representable   if and only if the defining  relations are unshrinkable words. In this section we will introduce a much more general  class of $*$-algebras which is defined  by imposing some conditions on the set of defining relations (see Definition~\ref{nonexpand}). For this class we will  prove $O\sp*$-representability. We also show that several unrelated, at first glance,  classes of $*$-algebras fall in this class.

 Denote by  $F_*$ a
  free associative algebra with generators $x_1,\ldots,
x_m$, $x^*_1,\ldots, x^*_m$. We do not incorporate the number of generators in the notations explicitly since it will be always clear from the context.   Algebra  $F_*$ is a $*$-algebra
with involution given on generators by $(x_j)^*=x_j^*$ for all $j
=1,\ldots, m$. Forgetting about involution we get  a free  associative algebra with $2m$ generators $F_{2m}$. The algebra  $F_*$ is
a semigroup algebra of a semigroup $W$ of all words in generators
$x_1,x_2,\ldots, x_m$, $x^*_1,x^*_2,\ldots, x^*_m$.

We have compiled all necessary prerequisites from  Gr\"obner basis theory of non-commutative  associative algebras in the appendix. Below we will explain how this theory will be applied for $*$-algebras.

 A set $S\subseteq
F$ of defining relations of an associative algebra $A$ is called
a Gr\"obner basis if it is closed under compositions (see Appendix).
 A \textit{Gr\"obner basis
of a $*$-algebra} $A$ is a Gr\"obner basis of $A$ considered as an
associative algebra. We need to put some extra requirements on a
Gr\"obner basis to make it "compatible" with the involution. The main
requirement we impose is a generalization of the notion of
unshrinkability of the word (see Definition~\ref{nonexpand} below).
\begin{dfn} A set $S\subseteq F_*$ is called
\textit{symmetric} if the ideal $\mathcal{I}$ generated by $S$ in
$F_*$ is a $*$-subalgebra of $F_*$.
\end{dfn} In particular, $S$ is
symmetric if $S^*=S$.

 For the notations $u\prec w$, ${ \rm R_S}(w)$, $BW$ and
order on $W$ used below we refer the reader to the appendix.
\begin{dfn}\label{nonexpand}
A symmetric subset  $S\subseteq F_*$ closed under compositions is called \textit{non-expanding} if for every $u, v, w \in BW $
such that $u\not= v$ and $ww^*\prec {\rm R}_{S}(uv^*)$ the
 inequality $w< \sup{(u,v)}$ holds, i.e. $w<u$ or $w<v$. If in addition for
every word $d \in BW$ the word $dd^*$ also belongs to  $BW$ then $S$ is called {strictly non-expanding}.
\end{dfn}
A $*$-algebra $A$ is called (strictly) non-expanding
 if it possesses a Gr\"obner basis $GB$ which is
 (strictly) non-expanding.
\begin{lem}
A symmetric closed under compositions subset $S\subseteq F_*$ is non-expanding if and only if for every $u, v\in BW$ such that $u>v$ and $\abs{u} = \abs{v}$ the property $uu^*\prec {\rm R}_{S}(uv^*)$ does not hold.
\end{lem}
\begin{proof}
Let for some $u, v, w\in BW$, $ww^*\prec {\rm R}_{S}(uv^*)$. Then  $ww^*\le uv^*$ and therefore  $\abs{w}\le \frac{\abs{u}+\abs{v}}{2}$. If $\abs{u}\not=\abs{v}$ then $\abs{w}< \max (\abs{u},\abs{v})$ and, consequently, $w<\sup(u,v)$. We can assume, henceforth, that  $\abs{u} = \abs{v}$.  Then $ww^*\le uv^*$ implies that $w\le u$. If $u<v$ then, clearly, $w<v$. If $u>v$ then by the assumptions of the Lemma $uu^*\not\prec {\rm R}_{S}(uv^*)$ and, hance, $w<u$ which proves the lemma.
\end{proof}

 Let $G \subseteq W$ and $T=[1,n]\cap\mathbb{Z}$ is an interval of positive integers with
$n=\abs{G}$. In case $\abs{G} = \infty$ we denote by $T$ the set of positive integers.  An
 {\it enumeration} of $G$  is  a bijection  $\phi:G \rightarrow
 T$   such that $u > v$ implies $\phi(u)>\phi(v)$. It is obvious that
enumerations exist for any given $G$.

Let  ${\rm H}:F_*\to F_*$ be a linear operator defined  by the
rule ${\rm H}(uu^*)=u$ for $u\in W$ and ${\rm H}(v)=0$ if $v$ is
not of the form $uu^*$ for some word $u$.

Fix  a set $S\subseteq F_*$ closed under compositions, an enumeration
$\phi:BW\to\mathbb{N}$ of the corresponding linear basis and a
sequence of positive real numbers $\xi=\{a_k\}_{k \in \mathbb{N}}$.
 Define a linear functional
$T_{\xi}^{\phi}:K\to \mathbb{C}$ by putting
$T_{\xi}^{\phi}(u)=a_{\phi(u)}$ for every word $u\in BW$, where $K$
denotes  the linear span of $BW$. Let $n=|BW|$  which can be infinite  and  $V$ denote  a vector space over $\mathbb{C}$ with a basis $\{e_k\}_{k = 1}^n$.

Define  $\sca{\cdot}{\cdot}_\xi$ to be a sesquilinear form on $V$ defined  by the following rules $$\sca{e_i}{e_i}_\xi = a_i,$$ and  $$\langle e_i,e_j \rangle_\xi = T_{\xi}^{\phi}\circ
{\rm H} \circ {\rm R}_{S}(uv^*),$$  where $\phi(u)=i$, $\phi(v)=j$, $u,v\in BW$.
   The definition is correct since $u$ and $v$  as above  are unique.

\begin{thm}\label{scal}
 If $S$ is strictly non-expanding then there exists a sequence
$\xi=\{ a_k \}_{k \in \mathbb{N}} \subset \mathbb{N}$ such that
the sesquilinear form $\langle \cdot , \cdot \rangle_{\xi}$ is
positively defined.
\end{thm}
\begin{proof}
  Let $g_{ij}=\langle e_i,e_j
\rangle_\xi$ for $i,\ j\in \mathbb{N}$ and let $G=(g_{ij})_{1\le
i,j\le\infty}$ denote the  Gram matrix.
We will use Silvester's criterion to show,  by induction on $m$, that $a_m$ can be chosen such that principal minor $\Delta_m
 > 0$. For $m=1$ put $a_1=1$ then $\Delta_1 = 1 > 0$. Assume that $a_1,\ldots,a_{m-1}$ are chosen such that $\Delta_1>0,\ldots, \Delta_{m-1} > 0$.

 By definition if $u \in BW$, then $uu^*$ is  also
in $BW$. Thus by  definition  we have $\langle
e_{\phi(u)},e_{\phi(u)} \rangle_\xi = a_{\phi(u)}$. Take some
$i\le m$ and $j\le m$  with $i\not=j$ and find unique $u, v\in BW$
such that $i=\phi(u),\ j=\phi(v)$. Then ${\rm R}_{S}(uv^*)=\sum_{k}\alpha_k
w_k$ for unique $\alpha_k\in \mathbb{C}$ and $w_k\in BW$. Clearly
$\langle e_{\phi(u)},e_{\phi(v)} \rangle_\xi$ is  $\sum_k
\alpha_k a_{\phi(h_k)}$ where the sum is taken over those $k$ for
which $w_k$ is of the form $w_k=h_kh_k^*$ for some word $h_k$. Since $S$ is non-expanding  we
have that $h_k<\sup{(u,v)}$. Hence $g_{ij}$ is a linear form in
variables $a_1,\ldots,a_{m-1}$. Decomposing determinant $\Delta_m$ by the $m$-th row we obtain $\Delta_m=\Delta_{m-1}a_m +p_{m}
(a_1,\ldots,a_{m-1})$ for some polynomial  $p_{m} \in \mathbb{C}[a_1,\ldots,
a_{m-1}]$. Since $\Delta_{m-1}>0$ it is clear that
$a_m$ can be chosen such that $\Delta_m>0$. This
 completes the inductive proof.
\end{proof}

The space $K$ is obviously isomorphic to $V$ via the map $u\to e_{\phi(u)}$.
Thus the inner product $\sca{\cdot}{\cdot}_\xi$ on $V$ gives rise
to an inner product on $K$ which will be denoted by the same
symbol. It is a routine to check that
$\sca{u}{v}_\xi=\alpha(P(u\diamond v^*))$, where $P:F_*\to F_*$  is the
projection on the linear span of positive words $W_+=W\cap F_{*+}$,  $\alpha:
K\to \mathbb{C}$ is a linear functional and  $\diamond$ is the operation defined in the appendix. Let $z\mapsto L_z$ denote the
right regular representation of $A=F_*/\mathcal{I}$, i.e.  $L_z(f)=f z$
for any $z$, $f \in A$.
\begin{thm}\label{main}
Let $S\subseteq F_*$ be strictly non-expanding and let $\mathcal{I}$ be  the ideal generated by $S$ in $F_*$. Then the right regular
representation $L$ of the quotient $*$-algebra $A=F_*/\mathcal{I}$ on a
pre-Hilbert space $(K,\sca{\cdot}{\cdot}_\xi)$ is a faithful
$*$-representation.
\end{thm}
\begin{proof}
The representation stated in the theorem is associated by the GNS
construction with the positive functional $\alpha(P(\cdot))$ on
$A$. Thus it is a $*$-representation. Indeed, as in the  GNS
construction the set $N=\{a \in A | \alpha(P(aa^*))=0\}$ is a right ideal
in $A$. We can define an inner product on $A/N$ by the usual rule
$\sca{a+N}{b+N}=\alpha(P(a^*b))$. It is easy to verify that the right
multiplication operators  define a $*$-representation of $A$ on pre-Hilbert
$A/N$. The only difference with classical  GNS construction is
that this representation could not be, in general, extended to the
completion of $A/N$.

We will show that this representation is  faithful. Take any $f={\sum_{i=1}^n c_i w_i}\in
A$, where $c_i \in \mathbb{C}, w_i \in BW$. Without loss of
generality consider $w_1$ to be the greatest word among $w_j$.
Then $L_f(w_1^*)$ contains element $w_1w_1^*$ with coefficient
$c_1$. Hence $L_f\not=0$.
\end{proof}
As a straightforward corollary of Theorem~\ref{main} we obtain the following.
\begin{cor}
Every strictly non-expanding $*$-algebra is $O\sp*$-re\-pre\-sent\-able.
\end{cor}

\section{Sufficient conditions of strictly non-expend\-ability. Examples.}\label{suffic}

 In this section we will show that the class of strictly non-expanding $*$-algebras contains several known classes  of $*$-algebras. To accomplish this we introduce below several other classes of $*$-algebras (see Definition~\ref{approp}, Corollary~\ref{cor}, and Theorem~\ref{kir}) and prove that they are contained in the class of non-expanding $*$-algebras. The definition given below may look complicated but, in fact, it is much easier to verify its conditions than the conditions of non-expanding $*$-algebra.  A more thorough look reveals that the conditions of Definition~\ref{approp}  and in the theorems in this section are algorithmically verifiable. In the end of the section we will present some concrete examples.

We  call a subset $S\subseteq F$ {\it  reduced} if for every $s\in
S$ and any word $w\prec s$ no word $\hat{s'}$ with $s'\in S$ is
contained in $w$ as a subword. If $S$ is closed under compositions then $S$ being  reduced is equivalent to ${\rm R_S} (s) = s$ for every $s\in S$.  If the set $S$ is closed under
compositions then one can obtain reduced set $S'$ closed under
compositions
 generating the same ideal by replacing each $s\in S$ with
 ${\rm R}_S(s)$.

\begin{dfn}\label{approp}
A symmetric reduced  subset $S\subseteq F_*$ is called  \textit{strictly
appropriate} if it is closed under compositions and for every $s
\in S$ and every word $u\prec s$ such that $|u|=\deg(s)$ the
following conditions hold.
\begin{enumerate}
\item   The word $u$ is unshrinkable. \item  If  $u\not=\hat{s}$,
$\hat{s}=a b$, and $u=a c$ for some words $b,c$ and nonempty word
$a$ then for any $s_1\in S$ such that there is word
$w\prec s_1,\ w\not= \hat{s}_1, |w|=|\hat{s}_1|$  either word
$\hat{s}_1$ does not contain $u$ as a  subword  or $\hat{s}_1$ and
$u$ do not form a composition in such a way that $\hat{s}_1=d_1 a
d_2$ and $u=a d_2 d_3$ with some nonempty words $d_1, d_2, d_3$.
\end{enumerate}
 A $*$-algebra $A$ is called  \textit{strictly appropriate} if it possesses a strictly appropriate Gr\"obner
basis.
\end{dfn}
We will use the following simple combinatorial facts proved in~\cite[Lemma~2]{pop1}.
For every two words $u$ and $v$ in free $*$-semigroup $W$ such that $uv^*=ww^*$ for some word $w$ either
$u=v$ or $v=udd^*$ for some $d\in W$ or $u=vcc^*$ for some $c\in W$ depending on  whether  $\abs{u}=\abs{v}$ or $\abs{u}<\abs{v}$ or  $\abs{u}>\abs{v}$.

  If $S=S^*$ is a closed under composition subset of $F_*$ such that $\hat{s}$ is unshrinkable for every $s\in S$ then $u\in BW$ if and only if $uu^*\in BW$.

 In the following theorem  for a word $w\in W$ of
even length $w=w_1w_2, |w_1|=|w_2|$ we will denote ${\rm H}_0(w)=w_1$.
\begin{thm}\label{principal}
Every strictly appropriate set $S\subseteq F_*$ is non-expanding. If
in addition  $S=S^*$ then $S$ is strictly non-expanding.
\end{thm}
\begin{proof}
Let $u, v\in BW$ be such that $u>v$ and $\abs{u} = \abs{v}$.

1. If $uv^*\in BW$ then $uu^*\prec {\rm R}_{S}(uv^*)$ implies
$uu^*=uv^*$ and, hence, $u=v$ which is a contradiction.

2. Now let $uv^*\not\in BW$. There are words $p,q \in BW$ and
element $s\in S$ such that $uv^*=p\hat{s}q$. Moreover, since $u,v
\in BW$ none of them can contain  $\hat{s}$ as a subword. Hence
$\hat{s}=a b$ with nonempty words $a$ and $b$ such that $u=p a$
and $v^*= b q$. Write down $s=\alpha \hat{s} + \sum_{i=1}^k \alpha_k w_i
+f$, where $w_i\in W$, $\alpha, \alpha_i\in\mathbb{C}$, and  $\deg(f) < \deg(s)$ and $|\hat{s}|=|w_i|$ for all  $i\in
\{1,\ldots, k\}$. Assume that for some integer $i$ word $pw_iq$
belongs to $BW$ and $pw_iq=uu^*$. If the middle
of the word $pw_iq$ comes across $w_i$, i.e.
$\max(|p|,|q|)<|u|$, then $w_i=c d$,  $u= p c$,  and $w^*=d q$ with
some nonempty words $c,d$. Hence $p c= q^*d^*$. If $|c|\le |d|$
then $d^*=g c$ for some word $g$ and so $w_i= c d = c c^* g^*$
which contradicts unshrinkability of $w_i$. If $|c|>|d|$ then $p
c= q^*d^*$ implies $c= g d^*$ for some word $g$ and we again see
that $w_i= g d^*d$ is shrinkable. Thus $\max(|p|,|q|)\ge |u|$. If
$|p|>|u|$ then $|u|=|p|+|a|>|u|$ which is impossible, hence   $|v|=|b|+|q|>|u|$.

3. Let $uv^*=p\hat{s}q$ and $s=\alpha \hat{s} + \sum_{i} \alpha_i w_i +f$
as above and $uu^*\prec {\rm R}_{S}(pw_iq)$ for some $i$.  Since $uu^*<p w_i q<uv^*$ word  $p w_i q$ begins
with $u$. If $\hat{s}=a b$ such that $p a=u, b q =v^*$ then $w_i$
 begins with $a$. Therefore  $\hat{s}$ and $w_i$ begin
with the same generator. Since $pw_i q\not\in BW$ there is
$s_1=\alpha_1 \hat{s}_1 +\sum_j\beta_j u_j +g\in S$ where $u_i\in W$, $\alpha_1, \beta_i\in\mathbb{C}$, and
$\deg(g)<\deg(s_1)$ such that $pw_iq=p_1\hat{s}_1q_1$ for some words
$p_1, q_1$. If we assume that  for some $j$ word $uu^*\prec {\rm
R}_{S}(p_1u_jq_1)$ then  ${\rm
H}_0(p_1u_jq_1)=u$ since $p_1u_jq_1<uv^*$. The word $\hat{s}_1$ can  not be a subword in the
first half of the word $pw_i q$ since ${\rm H}_0(p_1u_j q_1)={\rm
H}_0(pw_iq)=u$ and assuming the contrary we see that $\hat{s}_1$
and $u_j$ are both subwords of $u$ in the same position, hence they
must  be equal $\hat{s}_1=u_j$.
 The word $\hat{s}_1$ can not contain subword $w_i$ because of condition $2$ in the definition of strictly appropriateness.
 Obviously, $\hat{s}_1$ can not be a subword in $q$ because  $q\in
BW$. Thus either $w_i$ and $\hat{s}_1$  intersect (in the specified order) or $\hat{s}_1$
and  $w_i$ intersect in such a way that $\hat{s}_1=d_1 a d_2$ and
$w_i = a d_2 d_3$. But this contradicts the strictly
appropriateness of $S$. So we have proved that $S$ is non-expanding.
The fact that for any word $g\in BW$ word $gg^*$ lies in $BW$
follows from the remark preceding the theorem (see also~\cite[lemma~2]{pop1}).
\end{proof}

The following is a  simplification of the preceding theorem which is easier to verify in examples.
\begin{cor}\label{cor}
Let   $S\subseteq F_*$ be  symmetric and  closed under
compositions. Suppose that for every $s \in S$ and every word $u\prec s$ such that $|u|=\deg(s)$ the word $u$ is unshrinkable. In case $u\not= \hat{s}$ suppose also that
  words $\hat{s}$ and $u$  start with different generators. Then $S$ is non-expanding.
  If in addition $S=S^*$ then $S$ is strictly non-expanding.
\end{cor}

{\bf Example.}  Let $\mathcal{L}$ be a finite dimensional real Lie
algebra with  linear basis $\{e_j\}_{j=1}^{n}$. Then its universal
enveloping algebra $U(\mathcal{L})$  is a
$*$-algebra with involution given on generators as $e_j^*=-e_j$.
We claim that this $*$-algebra is non-expanding.
  Indeed $M=\{ e_i
e_j-e_j e_i- [e_i,e_j], i>j \}$ is a set of defining relations for
$U(\mathcal{L})$. It is closed under compositions (see example
in~\cite{Bok} or use PBW theorem).  Thus the set $S= \{e^*_j+e_j,
1\le j \le n \}\cup M$ is also closed under compositions (we
consider $e^*_1>e^*_2>\ldots>e^*_1>e_1>\ldots > e_n$) since
$e^*_j$ and $e_k e_l$ do not intersect for any $j$, $k$, $l$. It
is easy to see that $S$ is symmetric. Thus $S$ is non-expanding by
Corollary~\ref{cor}. However, $S\not=S^*$ and $S$ is not strictly non-expanding.

\begin{thm}\label{kir}
Let $S\subseteq F_*$ be a symmetric closed under compositions reduced  subset such that the following conditions are satisfied.
\begin{enumerate}
\item For every  $s\in S$  every word $w\prec s$ with
$|w|=\deg(s)$ is unshrinkable. \item For every $s_1$, $s_2\in S$
and every word $u\prec s_1$ with $|u|=\deg(s_1)$ the words $u$ and
$\hat{s_2}$ do not form a composition.
\end{enumerate}
Then  $S$ is non-expanding. If in
addition $S=S^*$ then $S$ is strictly non-expanding.
\end{thm}
\begin{proof}
Consider $u, v\in BW$ such that $u>v$ and $\abs{u} = \abs{v}$. We will prove that   $uu^*\not\prec {\rm R}_S(uv^*)$. Assume the contrary. Then
there is a sequence of words $\{q_i\}_{i=1}^{n}$ such that
$q_1=uv^*$, $q_n=uu^*$ and for every $1\le i\le n-1$ there is
$s_i\in S$ and words $c_i$, $d_i$, $u_i\in W$ such that $u_i\prec
s_i$, $u_i\not=\hat{s_i}$, $|u_i|=|\hat{s}_i|$ and $q_i=c_i
\hat{s}_i d_i$, $q_{i+1}=c_i u_i d_i$.

Let $j$ be the greatest with  the property that $\hat{s}_j$
intersects the middle of $q_j$. Such an index $j$ exists because
$j=1$ satisfies this property and we are making our choice  within a finite set.
Clearly $j<n$ since otherwise $u_{n-1}$ would be a subword in
$uu^*$ intersecting its middle and thus would be shrinkable, which contradicts assumption 1 of the theorem. Thus
for every $i\in \{j+1,\ldots,n-1 \}$ word $\hat{s}_i$ does not
intersect the middle of the word $c_{i-1} u_{i-1} d_{i-1}$. But
$\hat{s}_i$ could not be situated in the first half of this word
because otherwise the first half of the word $q_{i}$ would be
strictly less than $u$ and, consequently,  $q_n<uu^*$ which is a
contradiction. Thus $\hat{s}_i$ is a subword in the right half of
the word $q_{i}$. If $u_j$ and $\hat{s}_i$ does not form a composition for every $i\in \{j+1,\ldots,n-1 \}$ then $u_j$ is a subword in $uu^*$ intersecting its  middle and, thus, shrinkable. This contradicts  assumption~1 of the theorem. Hence $u_j$ and $\hat{s}_k$ intersect for some $k\in \{j+1,\ldots,n-1 \}$ contrary to assumption~2 of the theorem.  This proves that  $uu^*\not\prec {\rm R}_S(uv^*)$ and finishes the proof of the theorem.
\end{proof}

{\bf Examples}.

1. Let $S= \{w_j\}_{j\in \Re}$ be a set consisting of
unshrinkable words such that $S=S^*$. Since compositions of any two words are always
zero  this set is closed under compositions. The other conditions
in the definition of strictly non-expanding set is obvious. Thus
$*$-algebra
\[
\mathbb{C}\big\langle x_1,\ldots, x_n, x_1^*,\ldots,x_n^* | w_j,\
j\in \Re \big\rangle
\]
 is $O\sp*$-representable.

2.  Consider in more detail the simplest example of
monomial $*$-algebras ${A}_{x^2}=\algebra{x,x^*|x^2=0,
x^{*2}=0}$.

It was proved in~\cite{tapper} that $*$-algebra
$\algebra{x,x^*|x^p=0, x^{*p}=0}$ is $C^*$-representable for every
integer $p\ge 1$.  We will show  that among the representations of ${A}_{x^2}$ given
by Theorem~\ref{main} there is a $*$-representation in bounded
operators. It is an open problem for arbitrary $A_w=\algebra{x,x^*|w=0, w^{*}=0}$ with unshrinkable word $w$.

It can be easily verified that $BW$ consists of the words
$u_k=x(x^*x)^k$, $v_k=x^*(xx^*)^k$, $a_m=(xx^*)^m$, $b_m=(x^*x)^m$
where $k\ge 0, m\ge 1$. Obviously $BW_+$ consists  of the  words
$a_m$ and  $b_m$  ($m\ge 1$). If $z$, $w\in BW$ then $zw^*\in W_+$
if and only if $z$ and $w$ belong simultaneously to one of the sets $\{a_k\}_{k\ge 1}$, $\{b_k\}_{k\ge 1}$, $\{u_k\}_{k\ge 0}$, $\{v_k\}_{k\ge 0}$.
 Moreover,
$$ u_ku^*_t=a_{k+t+1},
v_kv^*_t=b_{k+t+1}, a_m a^*_n=a_{n+m},b_m b^*_n=b_{n+m}.$$
Consider the following ordering $$u_0<u_1<\ldots < a_1<a_2<\ldots
<v_0<v_1<\ldots < b_1<b_2<\ldots .$$ Denote
$\alpha(a_m)=\alpha_m$, $\alpha(b_m)=\beta_m$ then the Gram matrix
of the inner product defined in Theorem~\ref{scal} is  ${\rm diag}(A, A',B, B')$ where $A$, $A'$, $B$, $B'$ are Hankel matrices  $A=(\alpha_{i+j-1})_{ij}$,
$A'=(\alpha_{i+j})_{ij}$,  $B=(\beta_{i+j-1})_{ij}$,
 $B'=(\beta_{i+j})_{ij}$. Note that $Y'$ obtained from $Y$ by
canceling out the first column (here $Y$ stands for $A$ or $B$).

Thus the question of positivity of the form $\sca{\cdot}{\cdot}$
is reduced to the question of simultaneous positivity of two
Hankel matrices $C$ and $C'$ where the second is obtained from the first  by canceling out the first column. We will show that such matrices $A, A', B, B'$ could be chosen to be positive and such that $B=A$ and that the representation in Theorem~\ref{main} is in bounded
operators.

Let $f:[0,1]\to [0,1]$ be a continuous function such that $f(x)>0$ for all
$x\in [0,1]$. Let \[ \alpha_m = \int^1_0 t^{m+1} f(t) dt \] be the
moments of the measure with density $f(t)$.  It is well known that
the moment matrix $A=(\alpha_{i+j-1})_{i,j=1}^\infty$  is positively defined. But then $A'$ is the
moment matrix of the measure with density $t f(t)$ and thus is
also positively-defined. We can put $B=A$.

The representation acts on a Hilbert space $\rm H$ which is the completion of the linear space of the algebra $A_{x^2}$. Moreover, for all $k, t\ge 0$ and $m,n \ge 1$   \begin{align}\label{al} \sca{u_k}{u_t}= \alpha_{k+t+1}, \sca{v_k}{v_t}= \alpha_{k+t+1}, \\ \sca{a_m}{a_n}= \alpha_{m+n}, \sca{b_m}{b_n}= \alpha_{m+n+1}.\end{align}
All other inner products of the basis words are zero.

For every polynomial $P(t)= \sum_{k=0}^n c_k t^k$ with complex coefficients define $P(u)= \sum_{k=0}^n c_k u_k$ and, similarly, $P(v)= \sum_{k=0}^n c_k v_k$. If $c_0=0$ then we can define $P(a)= \sum_{k=0}^n c_k a_k$ and $P(b)= \sum_{k=0}^n c_k b_k$. Every element $g\in A_{x^2}$ can be expressed as $g= P(a)+Q(b)+R(v)+F(u)$ for some polynomials $P, Q, R, F$ such that $P(0)=Q(0)=0$.
To prove that the representation is in bounded operators we need
only to verify that the operator ${\rm L}_x$ of  multiplication  by  $x$ is  bounded. Obviously, $x u_k=0$ and $x a_m =0$  for all $k\ge 0$ and
$ m\ge 1$. Thus ${\rm L}_x(g)= Q(u)+R_1(a)$, where $R_1(t)=tR(t)$. Using (5)-(6) we obtain \begin{eqnarray*} \norm{{\rm L}_x(g)}^2 &=& \int_0^1 \abs{Q(t)}^2 t^2 f(t) dt + \int_0^1 \abs{R(t)}^2 t^3 f(t) dt,\\
\norm{g}^2 &=& \int_0^1 \abs{P(t)}^2 t f(t) dt + \int_0^1 \abs{Q(t)}^2 t f(t) dt +\\ && \int_0^1 \abs{R(t)}^2 t^2 f(t) dt+  \int_0^1 \abs{F(t)}^2 t^2 f(t) dt.
\end{eqnarray*}

  Thus $\norm{{\rm L}_x(g)}\le \norm{g}$. This proves that $L$ is a representation in bounded operators and, consequently, ${A}_{x^2}$ is $C^*$-representable.

3.  The $*$-algebra given by the generators and relations:
\[
\mathbb{C}\big<a_1,\ldots, a_n | a_i^*a_j=\sum_{k\neq l}
T_{ij}^{kl}a_la_k^*;i\neq j  \big>,
\]
 with  $T_{ij}^{kl}=\bar{T}_{ji}^{lk}$ is strictly non-expanding
by Corollary~\ref{cor}. Indeed, no two elements from defining relations form a composition and the greatest word of any relation
begins with some $a_j$ and all other words begin with some $a_k^*$.
 Hence this $*$-algebra is $O\sp*$-representable. Note that if  the additional relations
$a_i^*a_i=1 +\sum_{k,l} T_{ii}^{kl}a_la_k^*$ are imposed  we obtain  algebras allowing Wick ordering (see~\cite{jorgensen}).

4. Let  $S \subset \mathbb{C}W(x_1,\ldots,x_n)$  be closed under
 compositions  then a $*$-algebra

\[
A=\algebra{x_1,\ldots,x_n,x^*_1,\ldots,x^*_n|\  S\cup S^*}
\]
 is sometimes called $*$-double of $B= \algebra{x_1,\ldots,x_n|\  S}$.  By by
Corollary~\ref{dub} below $A$ is non-expanding.
 For finite dimensional algebra $B$ this already follows from Corollary~\ref{cor}. Indeed, if $S$ satisfies additionally the  property that  the greatest
word of every  relation begins with the generator different from the
beginnings of other longest words of this relation then $A$ is
strictly non-expanding by corollary~\ref{cor} since $S\cup S^*$ is, clearly,  closed under compositions. In particular, let $B$ be a  finite dimensional associative algebra with linear basis $\{e_k\}_{k=1}^n$. Then its "table of multiplication", i.e. the relations of the form  $e_i e_j - \sum c_{ij}^k e_k=0$, where $c_{ij}^k$ are the structure constants of the algebra $B$, forms a set of defining relations $S$ with the greatest words of length $2$ and others of
length $1$. Thus $*$-algebra $A=\algebra{x_1,\ldots,x_n,x^*_1,\ldots,x^*_n|\  S\cup S^*}$ is the  $*$-double of $B$. In other
words,  $A$ is the free product $B*B^*$,  where $B^*$ is an associative algebra such that there is an conjugate-linear anti-isomorphism $\phi:B\to B^*$, i.e. $\phi(a b) = \phi(b) \phi(a)$ and $\phi(\lambda a)= \overline{\lambda} a$ for $\lambda \in \mathbb{C}$ and $a, b\in B$. The involution on $A$ is given on the generators by the rules
$b^*=\phi(b)$ for any $b\in B$ and
$c^*=\phi^{-1}(c)$ for any $c\in B^*$.  The resulting $*$-algebra $A$ does not depend on the choice of $\phi$.

To deal with a general algebra  $B$ we need the following stronger result.
\begin{thm}\label{dubs}
Let $S=S^*$  be a closed
under compositions  subset of a free $*$-algebra $F_*$ with
generators $x_1,\ldots, x_n$, $x^*_1,\ldots, x^*_n$  such that for any $s\in S$ the following properties holds.
\begin{enumerate}
\item $\hat{s}\in G$ or
$\hat{s}\in G^*$ where $G=W(x_1,\ldots, x_n)$ is a semigroup generated by $x_1, \ldots, x_n$.
\item for any $u\prec s$ such that $|u|=|\hat{s}|$ words $u$ and
$\hat{s}$ both lie in the same semigroup $G$ or $G^*$.
\end{enumerate}
Then $S$ is strictly non-expanding.
\end{thm}
\begin{proof}
Let $X=\{x_1, \ldots, x_n\}$ and $X^*=\{x_1^*, \ldots, x_n^*\}$. As always $W$ will denote the semigroup $W(X\cup X^*)$.
If some word $w=y_{1}\ldots y_{t}$ where $y_{r}\in X\cup X^*$ contains subword $\hat{s}$ for some  $s\in S$ then  $w=p\hat{s}q$ for some words $p$ and $q$ in $W$. Let $ s= \hat{s} - \sum_{i=1}^n \alpha_i w_i$ ( $\alpha_i \in\mathbb{C}$, $w_i\in W$). The substitution rule
$\hat{s} \to \bar{s}$ (see the appendix)  replaces subword $w$ with $\sum_i \alpha_i p w_i q$. The conditions of the theorem ensure that all words $w_i$ such that $|w_i|=|\hat{s}|$ are in the same
semigroup either in $G$ or in $G^*$. Since decomposition ${\rm
R}_{S}(w)=\sum_j \beta_j u_j$, where $u_j\in BW$, $u_j=z^{(j)}_{1}\ldots
z^{(j)}_{k_j}$ with  $z^{(j)}_{r}\in X\cup X^*$ ($1\le r \le k_j$) can be obtained by several subsequent
substitutions considered above we see that for any $j$ such that
$|u_j|=|w|$ and for all $1\le r \le t$ both generators $z^{(j)}_{i_r}$ and $y_{k_r}$ are in the same set either $X$ or $X^*$.

Let $u,v\in BW$, $u>v$ and $\abs{u}=\abs{v}$. Assume that $ uu^*\prec
{\rm R}_{S}(uv^*)$. Without loss of generality we can assume that the word $u=z_1\ldots z_k$ ends with symbol from $X$, i.e. $z_k\in X$.
Then $uu^*=z_1\ldots z_k z_k^*\ldots z_1^* $.  By the first part of
the proof  $v^*$ begins
with a generator $x_l^*$ from the set $X^*$.
If $uv^*\not\in BW$ then there exists $s\in S$ such that $uv^* = p\hat{s} q$ for some words $p$ and $q$. Since $u, v\in BW$, $\hat{s}$ intersects both $u$ and $v^*$. Hence $\hat{s}$ contains $z_k x^*_l$  as a subword. This contradicts assumption~1 of the theorem.
 Thus $uv^*\in BW$ and ${\rm
R}_{S}(uv^*)=uv^*$. Clearly, $uv^*=uu^*$ implies $u=v$. Obtained
contradiction proves that $S$ is non-expanding.  Since for every $s\in S$, $\hat{s}$ is unshrinkable and    $S=S^*$  we have that  for any  $d\in BW$
word $dd^*$ is in  $BW$. Thus $S$ is strictly non-expanding.
\end{proof}

It could be shown using Zorn's lemma that for any algebra $A$ and
any its set of generators $X$ there is a Gr\"obner basis $S$ corresponding to $X$ with any given inductive ordering of the generators. It is easy to check that $S\cup S^*$ satisfies assumptions of Theorem~\ref{dubs}, thus,  we have the following.
\begin{cor}\label{dub}
If $B$ is a finitely generated  associative algebra then its
$*$-double  $A=B*B^*$ is strictly non-expanding $*$-algebra. Hence
$A$ has a faithful $*$-representation in pre-Hilbert space.
\end{cor}

Below we give some known examples of $*$-doubles which have finite
Gr\"obner bases.

5. We  present an example of $O\sp*$-algebra which is not $C^*$-re\-pre\-sent\-able.  Consider the  $*$-algebra:

\begin{align*} Q_{4,\alpha}=\mathbb{C} \big\langle
q_1,\ldots,q_4,q^*_1,\ldots , q^{*}_4|\ q_j^2=q_j, q_j^{*2}=q_j^*,\text{ for all $1\le j\le 4$},\\  \sum_{j=1}^4
q_j=\alpha,   \sum_{j=1}^4 q_j^{*}= \overline{\alpha }
\big\rangle.
\end{align*} which is the  *-double of the
algebra

\[
B_{n,\alpha}=\algebra{q_1,\ldots,q_4|\ q_j^2=q_j, \ \sum_j q_j=
\alpha }
\]
This algebra  has the following Gr\"obner basis:

$S=\{ q_1  q_1 - q_1, q_2 q_2 - q_2,
  -q_3 q_2 -2q_1 - 2q_2 - 2 q_3 + \alpha + 2 \alpha q_1+
      2 \alpha q_2 + 2  \alpha q_3 - \alpha^2 - q_1 q_2 - q_1 q_3 -
      q_2q_1 - q_2q_3 - q_3 q_1, q_3 q_3 -q_3,
  -q_3 q_1 q_2  -3\alpha + 5\alpha^2 - 2\alpha^3 +
      q_2(6 - 10\alpha + 4\alpha^2) +
      q_3(6 - 10\alpha + 4\alpha^2) +
      q_1(8 - 13\alpha + 5\alpha^2) + (3 -
            2\alpha)q_1  q_2 + (6 - 4\alpha) q_1 q_3 +
            (6 - 4\alpha)q_2 q_1 + (6 - 4\alpha) q_2 q_3 +
            (3 - 2\alpha)q_3 q_1 + q_1 q_2 q_1 +
      q_1 q_2q_3 + q_1q_3 q_1 + q_2q_1q_3 + q_2 q_3 q_1)
\} $.  More detailed treatment of this algebra can be found in ~\cite{sam,bart}.
Note that when $\alpha =0$ the $*$-algebra
$Q_{4,0}=B_{4,0}*B_{4,0}^*$ has only zero representation in bounded
operators (see~\cite{bart}). Thus for this $*$-algebra only representations in unbounded operators could exist.

The representability of $*$-algebras generated by projections connected by  linear relations is closely related to Horn type Spectral Problem \cite{SP1, SP2, SP3}. Such algebras have finite Gr\"obner
bases.

6. That the generators in the previous example are idempotents is
not important for $O\sp*$-representability,  we can consider the following example:
\begin{align*}T_{3,\alpha}=\mathbb{C} \big\langle q_1,q_2,q_3,q^*_1,q_2^* ,
q^{*}_3|\ q_j^3=q_j, q_j^{*3}=q_j^* \text{ for $1\le j\le 3$}, \\ \sum_j q_j=\alpha, \sum_j
q_j^{*}= \overline{\alpha } \big\rangle.
\end{align*}
It is the  $*$-double  of the algebra $\algebra{q_1,q_2,q_3|\
q_j^3=q_j, \ \sum_j q_j= \alpha }$. We will find its Gr\"obner basis. We have the following set of relations
$
 \{q_1^3 - q_1,
 q_2^3 - q_2,
 q_3^3 - q_3,
 q_1 + q_2 + q_3 - \alpha
 \}
$. From these relations it follows that this  algebra is generated by
$q_1$ and $q_2$. Thus we can consider the following equivalent  set of
relations: $ \{q_1^3 - q_1,
 q_2^3 - q_2,
(\alpha - q_1 - q_2)^3-(\alpha - q_1 - q_2) \}$.  Introduce the following  order on the generators $q_2>q_1$. All relations are
already  normalized, i.e. all leading coefficients are equal to $1$.  The greatest words in these relations are
$q_1^3$, $q_2^3$ and $q_1^2 q_2$. Thus we have no reductions. The
first and the third relations form  two compositions. From one side  they intersect by the word $q_1$. And the result of this
composition is
  $(q_1^3 - q_1)q_1q_2-q_1^2((\alpha - q_1 - q_2)^3-(\alpha - q_1 - q_2))$.
  On the other hand  they  intersect by the word $q_1^2$. The result of this composition is
  $(q_1^3 - q_1)q_2-q_1((\alpha - q_1 - q_2)^3-(\alpha - q_1 - q_2))$.
  Another composition is formed by the third and the second relations.
  Their greatest words intersect by the word $q_2$. Result of this  composition is
  $((\alpha - q_1 - q_2)^3-(\alpha - q_1 - q_2))q_2^2- q_1^2(q_2^3 -
  q_2)$.
Hence we have three new relations. After performing reductions we
will have the following set of relations:

$
S=\{ q_1^3 - q_1,
  -q_2^2 q_1 + 3\alpha q_1^2 + 3\alpha q_2^2 + \alpha^3 +
      q_1(-1 - 3\alpha^2) + q_2(-1 - 3\alpha^2) +
      3\alpha q_1  q_2 - q_1 q_2^2 - q_1^2  q_2 + 3\alpha q_2 q_1 -
      q_2 q_1^2 - q_1 q_2 q_1 - q_2 q_1 q_2, q_2^3 - q_2,
  -q_2 q_1  q_2  q_1^2 + -\alpha^3 + 9\alpha^5 -
      q_1^2(-3\alpha - 37\alpha^3) -
      q_2^2(3\alpha - 27\alpha^3) -
      q_2(-1 + 6\alpha^2 + 27\alpha^4) -
      q_1(18\alpha^2 + 30\alpha^4) - (-12\alpha -
            45\alpha^3)q_1 q_2 -
      27\alpha^2 q_1q_2^2 -( 1 + 30\alpha^2) q_1^2  q_2 +
      9\alpha q_1^2 q_2^2 - (6\alpha - 18 \alpha^3) q_2q_1
      - (1 + 3\alpha^2) q_2 q_1^2 -
      (-2 + 15\alpha^2) q_1q_2 q_1 +
      3\alpha q_1q_2 q_1^2 + 3\alpha q_1^2 q_2 q_1 -
      q_1^2q_2 q_1^2 - (-1 + 9\alpha^2) q_2 q_1  q_2 +
      6\alpha q_1q_2 q_1q_2 - q_1^2 q_2 q_1 q_2 -
      3\alpha q_2 q_1q_2q_1 + q_1q_2q_1q_2q_1  \}
$

Some of these relations do form compositions but all of them
reduce to zero. Hence it is a Gr\"obner basis. Thus $T_{3,\alpha}$ is $O\sp*$-representable for every complex parameter  $\alpha$.

\section{APPENDIX: Non-commutative Gr\"obner bases.}
For the convenience of the reader we review some relevant facts
from  non-commutative Gr\"obner bases theory (see~\cite{ufn,Bok})
with some straightforward  reformulations.

The reader should keep in mind that a Gr\"obner basis is just a special set of defining relations of a given algebra and thus is a subset of a free algebra. The main advantage of having a   Gr\"obner basis for an algebra is that one can algorithmically solve the  equality problem, i.e. one can decide for a given two non-commutative polynomial in the algebra generators if they represent the same element of the algebra or not.

The   Gr\"obner basis  always exists whatever system of generator one chooses but the  procedure to find a Gr\"obner basis does not always terminate.

Below we will present only those aspects of the Gr\"obner bases theory which are necessary for this paper.
 Let $W_n$ denote the  free semigroup  with   generators
 $x_1,\ldots, x_n$. For a word  $w= x_{i_1}^{\alpha_1} \ldots
 x_{i_k}^{\alpha_k}$ (where $i_1$, $i_2$, $\ldots$, $i_k\in \{1,\ldots, n\}$, and $\alpha_1$,
 $\ldots$, $\alpha_k \in \mathbb{N}\cup \{0\}$) the length of $w$,
 denoted by $|w|$, is defined to be $\alpha_1+\ldots +\alpha_k$.
  Let $F_n = \mathbb{C}\langle x_1, \ldots, x_n\rangle$ denote the free associative
algebra with   generators  $x_1,\ldots, x_n$. We will sometimes
omit subscript $n$.  Fix the  linear order on $W_n$ such that $x_1
> x_2 > \ldots > x_n$, the words of the same length ordered
lexicographically and the words of greater length are considered
greater. Any  $f \in F_n$  is a linear combination $\sum_{i=1}^k
\alpha_k w_i$  of distinct words $w_1$, $w_2$, $\ldots$, $w_k$ with complex coefficients $\alpha_i\not=0$ for all $i\in
\{1,\ldots, k \}$).  Let
$\hat{f}$ denote  the greatest of these words, say $w_j$. The coefficient $\alpha_j$  we denote by  ${\rm lc} (f)$ and call {\it leading coefficient}. Then
denote $\hat{f}- (\alpha_j)^{-1}f$ by $\bar{f}$. The degree of
$f\in F_n$, denoted by $\deg(f)$, is defined to be $|\hat{f}|$.
The elements of the free algebra $F$ can be identified with functions
$f:W\to \mathbb{C}$ with finite support via the map $f\to
\sum_{w\in W}f(w)w$. For a word $z\in W$ and an element $f\in F$
we will write $z\prec f$ if $f(z)\not= 0$.
\begin{dfn}
We will say that two elements $f,g \in F_n$ form a \textit{composition}
$w\in W$ if there are words $x,z \in W$ and nonempty word $y\in W$
 such that $\hat{f} = x y$, $ \hat{g} = y z$ and $w=x y z$. Denote the result of the composition $\beta f z-\alpha x g$ by
 $(f,g)_w$,  where $\alpha$ and $\beta$ are the leading coefficients of  $f$ and $g$ respectively.
\end{dfn}
If $f$ and $g$ are as in the preceding definition then  $f= \alpha x y +\alpha \bar{f}$ and $g= \beta y z +\beta \bar{g}$ and $(f,g)_w= \alpha \beta (\bar{f} z - x \bar{g})$. We will also say that $f$ and $g$ intersect by $y$. Note that there may exist many such $y$ for a given $f$ and $g$,  and the property "intersect" is not symmetrical. It is also obvious  that $(f,g)_w < w$. Notice  that two elements $f$ and $g$ may form compositions in many ways and $f$ may form composition with itself.

The following definition is due to
Bokut~\cite{Bok}.

\begin{dfn}
A subset $S \subseteq F_n$ is called \textit{closed under compositions}  if
for any two elements  $f$, $g \in S$ the following properties
holds.
\begin{enumerate}
\item If $f\not=g$ then the word $\hat{f}$ is not a subword in
$\hat{g}$. \item  If $f$ and $g$ form a composition $w$ then there
are words $a_j$, $b_j$ $\in W_n$,  elements
 $f_j \in S$  and complex   $\alpha_j$ such that $(f,g)_w=\sum_{j=1}^m \alpha_ja_jf_jb_j$ and  $a_jf_jb_j < w$, for $j=1,\ldots,m$.
\end{enumerate}
\end{dfn}

\begin{dfn}
A set $S\subseteq F$  is called a \textit{Gr\"obner basis} of an  ideal
$\mathcal{I}\subseteq F$ if for any $f\in \mathcal{I}$ there is
$s\in S$ such that $\hat{s}$ is a subword in $\hat{f}$. A
Gr\"obner basis $S$ of $\mathcal{I}$ is called minimal if no
proper subset of $S$  is a Gr\"obner basis of $\mathcal{I}$.
\end{dfn}

If $S$ is closed under compositions then $S$ is a minimal
Gr\"obner basis for the ideal $\mathcal{I}$ generated by $S$ (see
\cite{Bok}). Henceforth we will consider only minimal Gr\"obner
bases.  Thus we will say that $S$ is a Gr\"obner basis of an
associative algebra $A=F/\mathcal{I}$ if $S$ is closed under
composition and  generates  $\mathcal{I}$ as an ideal of $F$.
Let $GB$ be a Gr\"obner basis for $A$ and let  $\hat{GB}=\{ \hat{s} | s\in GB\}$. Denote by $BW(GB)$
the subset of those words in $W_n$ that contain no word from
$\hat{GB}$ as a subword. It is a well known fact that $BW(GB)$ is a
linear basis for $A$. Henceforth we will write simply $BW$ since
we will always deal with a fixed Gr\"obner basis.

If $S\subseteq F$ is closed under compositions and $\mathcal{I}$ is an  ideal generated by $S$ then each element $f+\mathcal{I}$ of the factor algebra $F/\mathcal{I}$ is the unique linear combination of basis vectors $\{w +\mathcal{I}\}_{w \in BW}$  $$ f+\mathcal{I} = \sum_{i=1}^n c_i (w_i+ \mathcal{I}). $$

We can define an operator ${\rm R}_S:F\to F$ by the following rule ${\rm R}_S(f) = \sum_{i=1}^n c_i w_i$. The element ${\rm R}_S(f)$ can be considered as a canonical form of the element $f$ in the factor algebra $F/\mathcal{I}$. Computing canonical forms we can algorithmically decide if two elements are equal in $F/\mathcal{I}$.

For example for a finite dimensional Lie algebra $\mathcal{L}$ with linear basis $\{e_i\}_{i\in M}$ and structure constants $C_{ij}^k$ ($[e_i,e_j] = \sum_{k} C_{ij}^k e_k$) the set of relations $e_i e_j - e_j e_i - \sum_{k} C_{ij}^k e_k$  with $i>j$ constitute a Gr\"obner basis for the universal enveloping associative algebra $U(\mathcal{L})$ and the canonical form is given by the PBW theorem.

Clearly  ${\rm R}_{S}$ is a
retraction on a subspace $K$ in $F$ spanned by $BW$. We can
consider a new operation on the space $K$:  $f\diamond g= {\rm
R}_{S}(f g)$ for $f$, $g \in K$. Then $(K,+,\diamond)$ becomes an
algebra which is isomorphic to $F/\mathcal{I}$.

Each element $s\in S$ in a Gr\"obner basis could be considered as a substitution rule $\hat{f} \to \bar{f}$ which tells us to replace each occurrence of the  subword $\hat{f}$ with $\bar{f}$. The canonical form ${\rm R}_{S}(f)$ can be computed step by step by performing all possible substitutions described above. The order in which the substitutions performed is not essential and only a finite number of substitutions could occur. From this it follows that if $w\prec R_{S}(u)$ for some words $w$ and
$u$ then $w<u$.
For example, take algebra $A=\mathbb{C}\langle a,\ b | ba = q ab \rangle$ for some complex $q$. Then considering $b>a$ we obtain that $S=\{ba - q ab\}$ is a Gr\"obner basis for $A$. We have only one substitution rule $ba \to q ab$. To obtain the canonical form of $b^2a$ we   compute $b (b a) \to q (b a) b \to q^2 b^2 a$. Thus ${\rm R}_{S}(b^2a) = q^2 b^2 a$.

\end{document}